\documentclass{amsart}
\usepackage{amsthm,amsmath,amsfonts}
\usepackage{mathrsfs}
\usepackage[colorlinks=false,linktocpage=true]{hyperref}
\usepackage{tocvsec2}
\usepackage{latexsym}
\newtheorem{thm}{Theorem}[section]
\newtheorem{cor}{Corollary}[section]
\newtheorem{lem}{Lemma}[section]
\newtheorem{prop}{Proposition}[section]

\theoremstyle{definition}
\newtheorem{defn}{Definition}[section]

\newtheorem{rem}{Remark}[section]
\newtheorem{notn}{Notation}[section]
\newcommand{\thmref}[1]{Theorem~\ref{#1}}

\newcommand{\lemref}[1]{Lemma~\ref{#1}}
\newcommand{\coref}[1]{Corollary~\ref{#1}}

\newcommand{\notnref}[1]{Notation~\ref{#1}}
\def\qed{\quad\vcenter{\hrule\hbox{\vrule height.6em\kern.6em\vrule}\hrule}}
\newenvironment{pf}{{\bigskip\textit{\newline Proof.}\quad}}{$\qed$\bigskip\newline}
\newenvironment{pf*}[1]{{\bigskip\textit{\newline#1.}\quad}}{$\qed$\bigskip\newline}
\numberwithin{equation}{section}
\font\tenscrpt=eusm10
\font\sevenscrpt=eusm10 scaled 700
\font\fivescrpt=eusm10 scaled 500
\newfam\eusmfam
\textfont\eusmfam=\tenscrpt
\scriptfont\eusmfam=\sevenscrpt
\scriptscriptfont\eusmfam=\fivescrpt

\def\s{\bold s}

\def\t{\bold t}

\def\Wt{W_{\t}}
\def\Ws{W_{\s}}

\def\Wpsx{W^{x}(\s)}
\def\Wx{W^{x}}

\def\Wxotnj{W^{x}_{\tnj}}

\def\Wxtnj{\Wxotnj(t_j)}

\def\vt{\bold t}
\def\p{\partial}

\def\nProdPartialut{\left(\frac{\p u}{ \p t_1}\right)\cdots\left(\frac{\p u}{\p t_{n}}\right)}

\def\eqdef{\overset{\triangle}{=}}

\def\R{\Bbb{R}}

\def\Nn{\Bbb{N}_{n}}
\def\Sn{\Bbb{S}_{n}}
\def\Rpn{{\Bbb{R}}_+^n}
\def\Rn{{\Bbb{R}}^n}
\def\bRpn{\partial{\Bbb{R}}_+^n}
\def\intRpn{\overset{\hspace{1mm}\circ\quad}{\Rpn}}

\def\Rpnmo{{\Bbb{R}}_+^{n-1}}
\def\Rp{{\Bbb{R}}_+}
\def\Rd{{\Bbb{R}}^{d}}
\def\P{\Bbb{P}}

\def\Px{\Bbb{P}_x}

\def\EP{{\Bbb E}_{\Bbb P}}

\def\EPx{\Bbb{E}_{\Px}}

\def\sR{\mathscr{R}}
\def\sE{\mathscr E}
\def\sU{\mathscr U}
\def\bU{\mathbb U}

\def\sF{{\mathscr F}}
\def\sFjt{{\mathscr F}^{(j)}_{t}}
\def\sFjsj{{\mathscr F}^{(j)}_{s_{j}}}

\def\OFP{(\Omega,\sF,\P)}

\def\sW{\mathscr W}
\def\sWj{\mathscr W^{j,x}}

\def\tnj{{{\bold t}_{\not j}}}

\def\snj{{{\bold s}_{\not j}}}

\def\utnj{u_{{\bold t}_{\not j}}}

\def\utx{u({\bold t},x)}

\def\Mtnsj{M_{\tnj}(s_j)}
\def\Mxtnsj{M^{x}_{\tnj}(s_j)}

\def\Wxtnsj{W^{x}_{\tnj}(s_j)}

\def\WxBnt{{\mathbb W}^{x}_{B^{(1)},\ldots,B^{(n)}}(\t)}
\def\WxBn{{\mathbb W}^{x}_{B^{(1)},\ldots,B^{(n)}}}
\def\dstyle{\displaystyle}
\def\ds{\displaystyle}

\def\i{\mathbf i}
\def\eqdef{:=}

\def\sm{\setminus}

\def\lap{\Delta}

\def\lbl#1{\label{#1}}

\def\pa{\partial}

\def\lqv{\left<}
\def\rqv{\right>}

\def\lab{\left|}
\def\rab{\right|}
\def\lpa{\left(}
\def\rpa{\right)}
\def\lbk{\left[}
\def\rbk{\right]}
\def\lbr{\left\{}
\def\rbr{\right\}}
\def\bdf{\begin{defn}}
\def\edf{\end{defn}}
\def\bcr{\begin{cor}}
\def\ecr{\end{cor}}
\def\bnt{\begin{notn}}
\def\ent{\end{notn}}
\def\brm{\begin{rem}}
\def\erm{\end{rem}}
\def\blm{\begin{lem}}
\def\elm{\end{lem}}
\def\bpf{\begin{pf}}
\def\bpfs{\begin{pf*}}
\def\epf{\end{pf}}
\def\epfs{\end{pf*}}

\def\beq{\begin{equation}}
\def\beqs{\begin{equation*}}
\def\eeq{\end{equation}}
\def\eeqs{\end{equation*}}
\def\bsp{\begin{split}}
\def\esp{\end{split}}
\def\bc{\begin{cases}}
\def\ec{\end{cases}}
\def\bt{\begin{tabular}}
\def\et{\end{tabular}}

\def\bthm{\begin{thm}}
\def\ethm{\end{thm}}
\def\bpr{\begin{prop}}
\def\epr{\end{prop}}
\def\babs{\begin{abstract}}
\def\eabs{\end{abstract}}
\def\ben{\begin{enumerate}}
\def\rencomalp{\renewcommand{\labelenumi}{(\alph{enumi})}}
\def\rencomrom{\renewcommand{\labelenumi}{(\roman{enumi})}}
\def\een{\end{enumerate}}
\def\KSstxy{K^{\textsc{KSS}(n,d)}_{\t;x,y}}

\def\KBstxy{K^{\textsc{BS}(n,d)}_{\t;x,y}}
\def\KBssxy{K^{\textsc{BS}(n,d)}_{\s;x,y}}
\def\KBmitzs{K^{\textsc{BM}}_{t_{i};0,s_{i}}}
\def\KBmjtzs{K^{\textsc{BM}}_{t_{j};0,s_{j}}}
\def\KBmjtzz{K^{\textsc{BM}}_{t_{j};0,0}}
\def\vsx{v({\bold s},x)}

\def\BT{{\mathbb X}^x_{B}}
\def\BTP{{\mathbb X}^x_{B}(t)}

\def\P{{\mathbb P}}
\def\E{{\mathbb E}}

\def\EP{{\mathbb E}_{\P}}

\def\R{{\mathbb R}}
\def\Rd{{\mathbb R}^d}

\def\Rp{{\mathbb R}_+}

\def\R{{\mathbb R}}

\def\BTP{{\mathbb X}^x_{B}(t)}
\def\sUj{\sU^{j}}

\def\sL{\mathscr L}
\begin{document}
\title[Brownian-time\hspace{-.3mm} Brownian\hspace{-.3mm} sheet\hspace{-.3mm} and\hspace{-.3mm} 4th\hspace{-.3mm} order\hspace{-.3mm} interacting\hspace{-.3mm} PDE{\scriptsize s}\hspace{-.3mm} systems]{From Brownian-time Brownian sheet to a fourth order  and a Kuramoto-Sivashinsky-variant interacting PDE{\scriptsize s} systems}
\author{Hassan Allouba}
\address{Department of Mathematical Sciences, Kent State University, Kent,
Ohio 44242}
\email{allouba@math.kent.edu}
\date{2/27/11}
\begin{abstract}  We introduce $n$-parameter $\Rd$-valued Brownian-time Brownian sheet (BTBS): a Brownian sheet where each ``time'' parameter is replaced with the modulus of an independent  Brownian motion.   We then connect BTBS to a new system of $n$ linear, fourth order, and interacting PDEs and to a corresponding fourth order interacting nonlinear PDE.  The coupling phenomenon is a result of  the interaction between the Brownian sheet, through its variance,  and the Brownian motions in the BTBS; and it leads to an intricate, intriguing, and random field generalization of our earlier Brownian-time-processes (BTPs) connection to fourth order linear PDEs.  Our BTBS does not belong to the classical theory of random fields; and to prove our new PDEs connections, we generalize our BTP approach in \cite{Abtp1,Abtp2} and we mix it with the Brownian sheet connection to a linear PDE system, which we also give along with its corresponding  nonlinear second order PDE and a $2n$-th order linear PDE that we also connect to Brownian sheet.  In addition, we introduce the $n$-parameter $d$-dimensional linear Kuramoto-Sivashinsky (KS) sheet kernel (or ``transition density''); and we link it to an intimately connected system of new linear Kuramoto-Sivashinsky-variant interacting PDEs, generalizing our earlier one parameter imaginary-Brownian-time-Brownian-angle kernel and its connection to a linear KS PDE.  The interactions here mean that our PDEs systems are to be solved for a family of functions, a feature shared with some well known fluids dynamics models.  The interacting PDEs connections established here open up another new fundamental front in the rapidly growing field of iterated-type processes and their connections to both new and important higher order PDEs and to some equivalent fractional Cauchy problems.   We connect the BTBS fourth order interacting PDEs system given here with an interacting fractional PDE system and further study it in another article.   
\end{abstract}
\maketitle
\newpage
\tableofcontents
\section{Introduction and results statements}
\subsection{A brief orientation and motivation}
In \cite{Abtp1,Abtp2}, we introduced a large class of interesting iterated processes that we called Brownian-time processes (BTPs)\footnote{A key of frequently used acronyms is provided in Appendix \ref{B} for the reader's convenience.}---stochastic processes in which the regular clock $t$ is replaced with the Brownian clock $|B(t)|$, where $B$ is a one-dimensional Brownian motion starting at $0$---and we connected them to novel fourth order PDEs.  The simplest such BTP is the Brownian-time Brownian motion (BTBM) $\BT\eqdef\lbr X^x(|B(t)|);0\le t<\infty\rbr$ where $X^x$ is an independent  Brownian motion starting at $x\in\Rd$; and we showed in \cite{Abtp1} that the associated PDE, for $\BT$ as well as for a large class of related BTPs, is
\begin{equation}\label{BTBMPDE}
\begin{cases}
  \dfrac{\partial}{\partial t}u(t,x)=\dfrac{\Delta_{x} f(x)}{\sqrt{8\pi t}}+\dfrac18\Delta_{x}^2u(t,x);& t>0,\,x\in\Rd,
\cr u(0,x)=f(x); & x\in\Rd.
\end{cases}
\end{equation}
I.e., if $u(t,x)=\E f(\BTP)$, then $u$ is the solution to the PDE \eqref{BTBMPDE} for an appropriate class of initial functions $f$ (see also  Allouba \cite{Abtp2}, DeBlassie \cite{DeB}, and Nane \cite{Nanesd}).  
The novelty in the PDE \eqref{BTBMPDE} is the ${\Delta_{x} f(x)}/{\sqrt{8\pi t}}$ term which expresses the BTP memory-preserving (non-Markovian) property.  Such PDEs connections are fascinating from many perspectives, here are some: 
\ben
\item probabilistically \eqref{BTBMPDE} is the PDE solved by ``running'' the intrinsically-interesting BTP, giving a dynamic description of the average of the modulated process $f(\BT)$; 
\item the term ${\Delta_{x} f(x)}/{\sqrt{8\pi t}}$ results in an interesting unconventional interaction between the initial data $f$ and the solution $u$ (this feature is under investigation for its potential applications to financial models); 
\item along with Zheng, we also initiated in \cite{Abtp1} the link between Brownian-time processes and fractional time derivatives by way of the BTP half-derivative generator---an implicit equivalence of \eqref{BTBMPDE} to temporally-fractional PDEs that was later made explicit and extended in several directions in the nice works of Nane,  Meerschaert, Baeumer, Vellaisamy, Orsingher, and  Beghin \cite{Nanesd,MNV09,BMN,BOap09,BOspa09}; and 
\item in \cite{Aks}, we adapted our BTPs construction in \cite{Abtp2,Abtp1} and their PDEs connections and gave an explicit representation of the solution to a linearization of the celebrated Kuramoto-Sivashinsky (KS) PDE in \textit{all spatial dimensions}.  Besides giving a probabilistic Brownian-time-motivated and explicit representation of this important equation, this achievement is noteworthy since the existence of the KS semigroup is still unsettled for dimensions larger than $1$ via traditional analytical methods, precluding the possibility to have definitive results using semigroup methods in higher dimensions.  We will have more to say about that in upcoming work. 
\een
In this article, we introduce a natural random field generalization of our Brownian-time Brownian motion in \cite{Abtp1,Abtp2}, which we call Brownian-time Brownian sheet (or BTBS); and we also introduce an interesting multi-time-parameter Kuramoto-Sivashinsky sheet kernel that generalizes our single-time parameter KS kernel in \cite{Aks,Abtpspde,Abtpsie}.  We then link them both to intriguing fourth order interacting PDEs.  These new PDEs connections add novel and intricate new features not present in their one parameter ($n=1$) $d$-dimensional ($d\ge1$) counterparts in \cite{Abtp1,Aks}.   Our interest in these fourth order systems of interacting PDEs is at least threefold (1) their connections to the intrinsically interesting BTBS random field in Subsection \ref{BTBSsec}; (2) their novel interaction terms which are fascinating, characteristic of the multi-parameter Brownian-time setting, and intricate nontrivial generalizations of the important one parameter case;  and (3) interacting PDEs are important in turbulence and fluids dynamics models as in the venerable second order Navier-Stokes equations in which the interaction
between pressure, velocity, and an external force (e.g., gravity) plays a crucial role.  In that light, we believe that the interactions in the PDEs of this article---which are probabilistically-naturally-derived multiparameter generalizations of our single-time
parameter, fourth order, and turbulent pattern formation equations---deserve more investigation.  
We state and further describe the two main results of this paper (\thmref{BTBSPDE} and \thmref{PDEKS}) in the following two subsections.  
\subsection{Brownian-time Brownian sheet and its fourth order interacting PDEs system}\lbl{BTBSsec}
An interesting and natural random field generalization of our BTBM in \cite{Abtp1,Abtp2} is the Brownian-time Brownian sheet, which we now introduce.  Let $B^{(1)},\ldots,B^{(n)}$ be $n$ independent copies of a standard one-dimensional Brownian motion  starting at $0$ and independent of an $n$-parameter ($n\ge1$) $\Rd$-valued Brownian sheet 
$$W^{0}=\lbr W^{0}(\t)=\lpa W^{0}_{1}(\t),\ldots,W^{0}_{d}(\t)\rpa; \t=(t_{1},\ldots,t_{n})\in\Rpn\rbr,$$ 
``starting'' at $0\in\Rd$ under $\P$---$\P\{W^{0}(\t)=0\}=1$ for $\t\in\bRpn$ (see \notnref{papernot})---all defined on a probability space $\OFP$.  Of course, the Brownian sheet coordinates $\lbr W^{0}_{1}(\t); t\ge0\rbr,\ldots,\lbr W^{0}_{d}(\t); t\ge0\rbr$ are assumed independent. For any $x=\lpa x_{1},\ldots,x_{d}\rpa\in\Rd$, let $$W^{x}:=W^{0}+x=\lbr W^{x}(\t)=\lpa W^{x_{1}}_{1}(\t),\ldots,W^{x_{d}}_{d}(\t)\rpa; \t=(t_{1},\ldots,t_{n})\in\Rpn\rbr.$$  We define the $n$-parameter $\Rd$-valued Brownian-time Brownian sheet (BTBS) based on $\Wx$ and $B^{(1)},\ldots, B^{(n)}$, and starting at $x\in\Rd$, by 
\beq\lbl{BTBSdef}
\bsp
\WxBnt\eqdef \Wx\lpa \lab B^{(1)}(t_{1})\rab,\ldots,\lab B^{(n)}(t_{n})\rab\rpa;\mbox{ } \t\in\Rpn.
\end{split}
\eeq
Clearly, $\P\lbk\WxBnt=x\rbk=1$ for $\t\in\bRpn$\footnote{We note that in \cite{Bu1} Burdzy discussed a process intimately related to the case $n=2$ of our BTBS in the context of a question he posed regarding the path properties of such a process.}.

Our first main result (\thmref{BTBSPDE} below) gives a novel system of linear fourth order interacting PDEs \eqref{btbsystem} and a corresponding nonlinear fourth order interacting PDE \eqref{btbsnonl}   that are solved by running the Brownian-time Brownian sheet random field  in \eqref{BTBSdef}.  When $n>1$, this new BTBS-PDE connection describes the dynamics of the average of the modulated Brownian-time Brownian sheet  $f\lpa\WxBn\rpa$ as it interacts with the averages of the intimately related random fields given by $ \lpa\prod_{i\in\Nn\setminus\lbr j\rbr}\lab B^{(i)}\rab\rpa^{2}f\lpa\WxBn\rpa$, for each $j\in\Nn:=\lbr1,\ldots,n\rbr$.   
Before stating our main result, it is helpful to adopt some simplifying notational conventions.
\begin{notn}\lbl{papernot} We always use the notation $\Nn:=\lbr1,\ldots,n\rbr$; and we denote by $\intRpn$ and $\bRpn$ the interior and boundary of $\Rpn$, respectively.  From this point on, we will alternate freely between $u_{\t}$, $u(\t)$, and $u(t_{1},\ldots,t_{n})$ for typesetting convenience and ease of exposition.  Moreover, $\tnj=(t_{i})_{i\in\Nn\setminus\{j\}}\in\Rpnmo$ will denote the $(n-1)$-tuple point in $\Rpnmo$ that is obtained from $\t=(t_{1},\ldots,t_{n})\in\Rpn$ by removing the $j$-th variable, $t_{j}$. The notations $u_{\t}$, $u(\t)$, and $\utnj(t_{j})$ will all mean the $n$-parameter function $u:\Rpn\to\Rd$ evaluated at $\t=(t_{1},\ldots,t_{n})$, and we use $\utnj(t_{j})$ whenever we need to highlight the variation in the variable $t_{j}$ for clarity.  The same comment applies for the notations $u(\t,x)$ and $u_{\tnj}(t_{j},x)$.
\end{notn}

We are now ready to state our first main result. 
\bthm[Brownian-time Brownian sheet: fourth order interacting PDEs connections]\lbl{BTBSPDE}
Let $\{\WxBnt; \t\in\Rpn\}$ be an $n$-parameter $\Rd$-valued BTBS based on a Brownian sheet $\Wx$ and Brownian motions $\lbr B^{(i)}\rbr_{i=1}^{n}$ and starting on $\bRpn$ at $x\in\Rd$ on $\OFP$.    
Let $f:\Rd\to\R$ be bounded and measurable such that all second order partial derivatives $\p^{2}_{kl}f$ are bounded and H\"older continuous, with exponent $0<\alpha\le1$, for $1\le k,l\le d$.   If 
\beq{}\label{BTBSexp0}
\bsp
u(\t,x)&=\EP\lbk f\lpa\WxBnt\rpa\rbk,
\\ \mathscr{U}^{(j)}(\t,x)&=\EP\lbk \lpa\prod_{i\in\Nn\setminus\lbr j\rbr}\lab B^{(i)}(t_{i})\rab\rpa^{2}f\lpa\WxBnt\rpa\rbk;\ j\in\Nn,
\end{split}
\eeq
for $(\t,x)\in\Rpn\times\Rd$, then the family $\lbr u(\t,x),  \sU^{(1)}(\t,x),\ldots, \sU^{(n)}(\t,x)\rbr$ is a solution to both
\ben
\rencomrom
\item the system of fourth order interacting Linear PDEs:
\begin{equation}\label{btbsystem}
\begin{cases}
(a)\ \displaystyle\frac{\p u}{\p t_j}=\sqrt{\frac{\prod_{i\in\Nn\setminus\lbr j\rbr}t_i}{2^{4-n}\pi^{n}t_{j}}}\lap_{x} f(x)+\frac18\lap_{x}^{2}\sU^{(j)};& \t\in\intRpn,x\in\Rd,\\
(b)\ u(\t,x)=f(x);&\t\in\bRpn,x\in\Rd,\\
(c)\ \sU^{(j)}(\t,x)=0;&\tnj\in\p\Rpnmo, x\in\Rd,\\
(d)\ \sU^{(j)}(\t,x)=\lbk \prod_{i\in\Nn\setminus\lbr j\rbr}t_{i}\rbk f(x);&\ t_{j}=0, x\in\Rd,
 \end{cases}
\end{equation}
for $ j\in\Nn$$;$ and
\item the corresponding fourth order nonlinear interacting PDE:
\begin{equation}\label{btbsnonl}
\bsp
\prod_{j\in\Nn}\frac{\p u}{\p t_{j}}=\ &\sqrt{\frac{(t_1\cdots t_{n})^{n-2}}{2^{4n-n^{2}}\pi^{n^{2}}}}\lpa\lap_{x}f(x)\rpa^{n}
+\frac1{8^{n}}\prod_{j\in\Nn}\lap_{x}^{2}\sU^{(j)}
\\&+\sum_{\lpa k_{1},\dots,k_{n}\rpa\in\Sn}\prod_{j\in\Nn}T_{k_{j},j}(\t,x);\ \t\in\intRpn, x\in\Rd
\end{split}
\end{equation}
together with the temporal boundary conditions $(b)$, $(c)$, and $(d)$ in \eqref{btbsystem},
where 
\beq\lbl{crossterms}
T_{k_{j},j}(\t,x)=\bc\ds\sqrt{\frac{\prod_{i\in\Nn\setminus\lbr j\rbr}t_i}{2^{4-n}\pi^{n}t_{j}}}\lap_{x} f(x);&k_{j}=1,j\in\Nn\\
\ds\frac18\lap_{x}^{2}\sU^{(j)};&k_{j}=2, j\in\Nn,
\ec
\eeq 
and where $$\Sn=\lbr\lpa k_{1},\ldots,k_{n}\rpa;k_{1},\ldots,k_{n}\in\{1,2\}, n+1\le\sum_{j=1}^{n}k_{j}\le 2n-1\rbr, \ n\ge2.$$
\een
\ethm

The intriguing interaction (or coupling) in \thmref{BTBSPDE} between $u$ and $\sU^{(j)}$, $j=1,\ldots,n$, shows that the PDEs in \eqref{btbsystem} and \eqref{btbsnonl} are nontrivial intricate generalizations of our BTP-PDE connections in \cite{Abtp1,Abtp2} (the case $n=1$).  This precise coupling phenomenon is caused by the interaction between the Brownian-times $\lab B^{(1)}\rab,\ldots,\lab B^{(n)}\rab$ and the variance of the Brownian sheet $\Wx$ through $\pa_{j}\mbox{Var}\lpa\Wx(t)\rpa=\prod_{i\in\Nn\sm\{j\}}t_{i}$.  The nonlinear  feature of the BTBS-PDE connection and its corresponding linear system connection are inherited from the nonlinear Brownian sheet PDE together with its corresponding linear system of PDEs, which we give below in \lemref{BSPDE}.  However, in addition to the new coupling above, and due to the iterative nature of BTBS, those PDEs in \eqref{btbsystem} and \eqref{btbsnonl} are fourth order versus the much simpler second order Brownian sheet nonlinear PDE and linear system of PDEs.      On the other hand, when $n=1$, our BTBS becomes a Brownian-time Brownian motion (see \cite{Abtp1,Abtp2}), $u=\sU^{(j)}$, and \eqref{btbsystem} and \eqref{btbsnonl} reduce to our fourth order BTP-PDE \eqref{BTBMPDE} ((0.2) in \cite{Abtp1}).
 Note that if $n=1$ then $\prod_{i\in\Nn\setminus\lbr j\rbr}t_i=\prod_{i\in\Nn\setminus\lbr j\rbr}\lab B^{(i)}(t_{i})\rab^{2}=1$ since $\Nn\sm\lbr j\rbr=\phi$ (the empty set), and $\sum_{\lpa k_{1},\dots,k_{n}\rpa\in\Sn}\prod_{j\in\Nn}T_{k_{j},j}(\t,x)=0$ since $\Sn=\phi$.  Moreover, the interaction in the \eqref{btbsystem} means that the PDEs are to be solved for a family of functions $\lbr u,\sUj\rbr_{j\in\Nn}$ as opposed to just one function $u$ as in either the Brownian sheet case \eqref{heatsystem} or the one-parameter case \eqref{BTBMPDE}.
  
 As with its one-parameter BTP counterpart, the BTBS is memory preserving as is indicated by the inclusion of the Laplacian of the initial function $\lap_{x}f(x)$ in the BTBS PDEs in \eqref{btbsystem} and \eqref{btbsnonl}.  Also, just as BTPs are not classical (not semimartingales, not Markovian, and not Gaussian), BTBS is not a classical random field. To prove \thmref{BTBSPDE} we generalize and mix our Brownian-time approach from \cite{Abtp1,Abtp2} with the Brownian-sheet-PDEs connection which we give in \lemref{BSPDE} below.

\thmref{BTBSPDE} gives intriguing connections between the Brownian-time Brownian sheet and a novel nonlinear fourth order and interacting PDE \eqref{btbsnonl} and its corresponding interacting system \eqref{btbsystem}.     This connection opens up another new fundamental front in the rapidly growing field of Brownian-time and iterated  processes and their connections to both (1) new and important higher order deterministic and stochastic PDEs and (2) some equivalent fractional Cauchy problems (e.g., Allouba and Allouba et al.~\cite{Abtp1,Abtp2,Aks,Abtpsie,Abtpspde,Aksspde,ALks,AX,ADsh} and followup papers; as well as  the recent articles by Nane \cite{Nanetr,Nanesd}, Baeumer et al.~\cite{BMN},  Meerschaert et al.~\cite{MNV09}, Orsingher et al.~\cite{BOap09,BOspa09}, and DeBlassie \cite{DeB}).  In this new Brownian-time Brownian sheet PDEs direction, the BTBS-PDEs connections \eqref{btbsnonl} and \eqref{btbsystem} introduce a way for a random field to shed light on nonlinear PDEs---and their corresponding system of $n$ linear PDEs---that are fourth order and interacting; and these PDEs, in turn,  could be used to gain insight into the underlying Brownian-time random field.  Guided by our one-time-parameter  ($n=1$) work in \cite{Abtp1,Abtp2,Aks,Aksspde,Abtpsie}, we explore ramifications of such a program in planned upcoming work. 
In addition, following our work with Zheng in \cite{Abtp1}, where we initiated the link between Brownian-time processes and fractional time derivatives by way of the BTP half-derivative generator---an implicit connection to temporally-fractional PDEs that was later made explicit and extended in several directions in the works of Nane, Meerschaert, Baeumer, Vellaisamy, Orsingher, and  Beghin \cite{Nanesd,MNV09,BMN,BOap09,BOspa09}---we intend to connect our BTBS with systems of fractional interacting PDEs  and their corresponding nonlinear fractional coupled PDE in future work. 
\subsubsection{Three Brownian sheet-PDEs links}
Here, we give the connections between the $n$-parameter $\Rd$-valued Brownian sheet $\Wx$ and (1) a system of $n$-linear second order PDEs, (2) a corresponding nonlinear second order PDE, and (3) a linear $2n$-th order PDE.  This Brownian sheet-PDEs connection, which characterizes Brownian sheet similar to the way the heat equation characterizes Brownian motion, is one of the needed ingredients in our proof of \thmref{BTBSPDE}; and it contains within it the first element causing the coupling between $u$ and $\sU^{j}$ in \thmref{BTBSPDE} above.  Although some of its aspects are perhaps known as ``folklore'' among the specialists (we obtained the Brownian sheet-PDEs connection---as part of our unpublished preprint \cite{Ahyp}---and we circulated it as early as 1996), this connection does not seem to be recorded explicitly anywhere in the standard literature.   We present it here, along with a proof that we give in the Appendix for completeness.  From a mathematical finance point of view, the system \eqref{heatsystem} below is at the heart of our new Black-Scholes-Merton system for semi-SPDEs bonds and stocks models that we introduce and study in a separate joint paper.    Due to the fundamental nature of Brownian sheet and its connection to space-time white noise, the Brownian sheet-PDE connection in \lemref{BSPDE} below (as well as \coref{bsnonso} and \coref{bslpdeho})---a connection which is as of yet really untapped---and natural generalization thereof are useful to further study the interactions between random fields and SPDEs solutions on one hand and nonlinear PDEs or higher order linear PDEs  and their corresponding system of linear PDEs on the other (more on that in \cite{Ahyp}).  

The following result gives the Brownian sheet system of PDEs.
\blm[Brownian sheet linear second order PDEs system]\label{BSPDE}
Let $\{W^{x}(\t); \t\in\Rpn\}$ be an $n$-parameter $\Rd$-valued Brownian sheet on $\OFP$ with $\P\{W^{x}(\t)=x\}=1$,
for $\t\in\bRpn$.  Let $u(\t,x)=\EP[f(\Wx(\t))]$, for a continuous bounded real-valued function
$f:\Rd\to\R$.  Then, $u(\t,x)$ is a bounded solution to  the system of heat-type PDEs
\begin{equation}\label{heatsystem}
\begin{cases}
(a)\ \ds\frac{\p u}{\p t_j}=\left(\frac12{\prod_{i\in\Nn\setminus\lbr j\rbr}}t_i\right)\lap_{x}u;& j\in\Nn, \t\in\intRpn, x\in\Rd,\\
(b)\ u(\t,x)=f(x);&\t\in\bRpn,  x\in\Rd,
\end{cases}
\end{equation}
Moreover, $u(\t,x)$ is the unique bounded solution to 
\begin{equation}\label{jthheatpde}
\begin{cases}
(a)\ \ds\frac{\p u}{\p t_j}=\left(\frac12{\prod_{i\in\Nn\setminus\lbr j\rbr}}t_i\right)\lap_{x}u;& \t\in\intRpn, x\in\Rd,\\
(b)\ u_{\tnj}(0,x)=f(x);&\tnj\in\Rpnmo,  x\in\Rd,
\end{cases}
\end{equation}
for each $j\in\Nn$.  In particular, $u(\t,x)$ is the unique bounded solution to the system \eqref{heatsystem}.
 \elm
 As immediate corollaries to \lemref{BSPDE} we connect the Brownian sheet $\Wx$ to both a nonlinear second order PDE and a higher order linear PDE.
 \bcr[The Brownian sheet second order nonlinear PDE]\lbl{bsnonso}
 Under the same conditions in \lemref{BSPDE} $u(\t,x)=\EP[f(\Wx(\t))]$ is a bounded solution to the nonlinear second order PDE
\begin{equation}\label{nonl}
\begin{cases}
\displaystyle\nProdPartialut=\frac{(t_1\cdots t_{n})^{n-1}}{2^n}\lpa\lap_{x}u\rpa^{n};
&\t\in\intRpn, x\in\Rd,\\
u(\t,x)=f(x);&\t\in\bRpn, x\in\Rd,
\end{cases}
\end{equation}
\ecr
\bcr[The Brownian sheet $2n$-th order linear PDE]\lbl{bslpdeho}
Under the same conditions in \lemref{BSPDE} $u(\t,x)=\EP[f(\Wx(\t))]$ is a bounded solution to the $2n$-th order linear PDE
\begin{equation}\label{holpde}
\begin{cases}
\displaystyle\frac{\partial^{n}u}{\p t_{1}\cdots\p t_{n}}=\sL_{n}(\t,u);
&\t\in\intRpn, x\in\Rd,\\
u(\t,x)=f(x);&\t\in\bRpn, x\in\Rd,
\end{cases}
\end{equation}
where, for each $\t$, the $2n$-th order operator $\sL_{n}(\t,\cdot)$, acting on $u$, is obtained by using the Brownian sheet PDE system \eqref{heatsystem} and iteratively applying the $t_{j}$ derivatives for $j\in\Nn$.  In particular, if $n=1,\ldots,4$, then
\beq
\bsp
\sL_{1}(\t,u)&=\frac12\Delta u\\
\sL_{2}(\t,u)&=\frac12\Delta u+\frac{t_{1}t_{2}}4\Delta^{2} u\\
\sL_{3}(\t,u)&=\frac12\Delta u+\frac34 t_{1}t_{2}t_{3}\Delta^{2} u+\frac{t_{1}^{2}t_{2}^{2}t_{3}^{2}}8\Delta^{3} u\\
\sL_{4}(\t,u)&=\frac12\Delta u+\frac84 t_{1}t_{2}t_{3}\Delta^{2} u+\frac{5}8t_{1}^{2}t_{2}^{2}t_{3}^{2}\Delta^{3} u+\frac{t_{1}^{3}t_{2}^{3}t_{3}^{3}}{16}\Delta^{4} u\\
\end{split}
\eeq
\ecr
  \brm
  \ben\rencomalp
  \item The uniqueness assertion in \lemref{BSPDE} is not needed in our proof of \thmref{BTBSPDE}, but we give a uniqueness proof in Appendix \ref{app} that judiciously exploits the relationship between Brownian sheets and their corresponding space-time white noises to essentially reduce the proof to stochastic analysis in one time parameter. 
\item The uniqueness for \eqref{jthheatpde} for each $j\in\Nn$ together with the fact that $u(\t,x)$ solves \eqref{heatsystem} means that, in the Brownian-sheet case, the system \eqref{heatsystem} is overdetermined or reducible.  I.e., reducing the system \eqref{heatsystem} to \eqref{jthheatpde} (for any $j$) leaves the bounded-solutions set unchanged.
 \item The Brownian sheet heat-type system \eqref{heatsystem} is symmetric in $\t=\lpa t_{1},\ldots,t_{n}\rpa$:  if $\pi(\t)$ is any one of the $n!$ permutations of $\t=\lpa t_{1},\ldots,t_{n}\rpa$, then $u(\t,x)$ and $u_{\pi}(\t,x)=u(\pi(\t),x)$ satisfy \eqref{jthheatpde} for each $j\in\Nn$, without change.  This is not surprising since the unique solution $u(\t,x)=\EP[f(\Wx(\t))]$ is symmetric in $\t$.  Restricting our attention to $\t$-symmetric bounded solutions $s(\t,x)$ to the nonlinear PDE \eqref{nonl}, we obtain that $s(\t,x)$ must satisfy \eqref{jthheatpde} for each $j$, and hence $s$ is the unique solution $s(\t,x)=\EP[f(\Wx(\t))]$.  A similar observation holds for the higher order  PDE \eqref{holpde}.  We leave the details to the interested reader.
 \een
 \erm
 Unlike the Brownian motion, ``running'' the Brownian sheet $\Wx$ solves heat PDEs with variable conductivity.  The system of heat-type PDEs \eqref{heatsystem} is one in which the $j$-th PDE (where  time is given by $t_j)$ may be thought of as modeling heat flow in a $d$-dimensional medium whose {\sl thermal diffusivity} (conductivity) is determined by the product of all other $n-1$ parameters $t_i$, $i\ne j$, in the driving Brownian sheet.  I.e., the conductivity of the $d$-dimensional medium varies with $n-1$ outside factors represented by the parameters $t_i$, $i\ne j$, in the driving Brownian sheet.  Probabilistically, this is a manifestation of the fact that $\pa_{j}\mbox{Var}(\Wt^{x})=\prod_{i\in\Nn\setminus\lbr j\rbr}t_i$---which interacts with the Brownian times $\lab B^{1}\rab,\ldots,\lab B^{n}\rab$ in the BTBS $\WxBnt$ to result in the coupling between $u$ and $\sU^{j}$ in \eqref{btbsystem} and \eqref{btbsnonl}.    Naturally, when $n=1$, the higher order PDE \eqref{holpde}, the nonlinear PDE \eqref{nonl}, and the corresponding system \eqref{heatsystem} reduce to the standard heat equation associated with Brownian motion.   

\subsection{A new intimately connected Kuramoto-Sivashinsky-variant  interacting PDEs system}
Generalizing our KS links to variants of BTPs in \cite{Aks,Aksspde,Abtpspde,Abtpsie}, we modify the BTBS setting above to obtain solutions to interesting systems of novel interacting variants of the celebrated KS equation.  For a sample of some deterministic KS results, background, and its applications to turbulence phenomena in chemistry and combustion the reader may consult \cite{T} and the references therein.  
In what follows $\i=\sqrt{-1}$ and should not be confused with the index variable $i$.
In our second main result (\thmref{PDEKS}); we show that an interaction similar to that obtained for the BTBS-PDEs  is also possible in the Kuramoto-Sivashinsky setting, we give the corresponding system of linear Kuramoto-Sivashinsky-variant interacting PDEs in \eqref{KSPDE}, and we give an explicit representations for the three-components-solution to each PDE in terms of the complex multiparameter KS kernel 
\beq\lbl{kssker}
\KSstxy:=\int_{\Rn}\exp\left(\i\sum_{i=1}^{n}s_{i}\right)p_{\i\s}^{(d)}(x,y)\prod_{i=1}^{n}\KBmitzs d\s
\eeq  
which we call the ``transition density'' (or kernel) of the $n$-parameter $d$-dimensional linear Kuramoto-Sivashinsky sheet (LKSS).  The kernel 
\beq\lbl{propag}
p_{\i\s}^{(d)}(x,y)=\frac{e^{-|x-y|^2/2\i\prod_{i=1}^{n}s_i}}{\lbk{2\pi \i \prod_{i=1}^{n}s_i}\rbk^{d/2}}; \ x,y\in\Rd,\ \s=(s_{1},\ldots,s_{n})\in\Rn.
\eeq
 is the $n$-parameter generalization of the $d$-dimensional famous propagator (complex Gaussian density on $\Rd$) in \cite{Aks}, associated with Schr\"odinger equation, and $\KBmitzs$ is the density of the BM $B^{i}$ for each $i\in\Nn$.     The kernel $\KSstxy$ naturally generalizes the important linear Kuramoto-Sivashinsky kernel in \cite{Aks,Aksspde,Abtpspde,Abtpsie} (the case $n=1$ for $d\ge1$) associated with the linear KS PDE.  Working with the complex variant of the BTBS density in \eqref{kssker}, we now briefly provide a link to such deterministic KS-variant interacting  system. 

Motivated by the BTBS case above and our KS link in \cite{Aks}, we define
\begin{equation}\label{expectn}
\begin{split}
v(\s,x)&\eqdef\exp\left(\i\sum_{i=1}^{n}s_{i}\right)\int_{\Rd}f(y)p_{\i\s}^{(d)}(x,y) dy\\
u(\t,x)&=\int_{\Rn} v(\s,x) \prod_{i=1}^{n}\KBmitzs d\s\\
\bU^{(j)}(\t,x)&=\int_{\Rn} \lpa\prod_{i\in\Nn\setminus\lbr j\rbr}s_i\rpa v(\s,x) \prod_{i=1}^{n}\KBmitzs d\s\\
\sU^{(j)}(\t,x)&=\int_{\Rn} \lpa\prod_{i\in\Nn\setminus\lbr j\rbr}s_i\rpa^{2}v(\s,x) \prod_{i=1}^{n}\KBmitzs d\s
\end{split}
\end{equation}
for $(\t,x)\in\intRpn\times\Rd$ and $j\in\Nn$.  
To define $u$, $\bU^{(j)}$, and $\sU^{(j)}$ on $\bRpn\times\Rd$ we adopt the following notational convenience: for $\t\in\Rpn$ and for $I\neq\phi$ and $I\subset\Nn$ we denote by $\zeta_{I}(\t)$ the point in $\Rpn$ obtained from $\t=(t_{1},\ldots,t_{n})$ by setting  $t_{k}=0$ for every $k\in I$.  For any function $g$ on $\Rpn\times\Rd$, $g_{\zeta_{I}}(\t,x):=g(\zeta_{I}(\t),x)$.  With this in mind, we give the following temporal boundary definitions:
\beq\lbl{ksbdrycnd}
\bsp
u_{\zeta_{I}}(\t,x)&=f(x)\exp\lpa{\ds-\frac12\sum_{k\in I^{c}}t_{k}}\rpa;\quad I\subset\Nn,\\
\bU^{(j)}_{\zeta_{I}}(\t,x)&=\bc f(x)\lbk\ds\prod_{k\in\Nn\setminus\{j\}}\i t_{k}\rbk \exp\lpa{\ds-\frac12\sum_{k\in\Nn\setminus\{j\} }t_{k}}\rpa;&I=\lbr j\rbr,\cr
0;&I\subset\Nn\setminus\lbr j\rbr,
\ec\\
\sU^{(j)}_{\zeta_{I}}(\t,x)&=\bc f(x)\lbk\ds\prod_{k\in\Nn\setminus\{j\}}\lpa t_{k}-t_{k}^{2}\rpa\rbk \exp\lpa{\ds-\frac12\sum_{k\in\Nn\setminus\{j\} }t_{k}}\rpa;&I=\lbr j\rbr,\cr
0;&I\subset\Nn\setminus\lbr j\rbr,
\ec
\end{split}
\eeq
for $I\neq\phi$.  Of course, if $I=\Nn$ in \eqref{ksbdrycnd} then $I^{c}=\phi$ and $u_{\zeta_{I}}(\t,x)=f(x)$.  Also, if $n=1$ then $\sU^{(j)}_{\zeta_{I}}(\t,x)=\bU^{(j)}_{\zeta_{I}}(\t,x)=u_{\zeta_{I}}(\t,x)=f(x)$.

We are now ready to state our second main result.
\begin{thm}[Linear KS-variant interacting PDEs system]\label{PDEKS}
Let $f\in C^2_c(\Rd;\R)$ such that the partial derivatives $\partial_{ij} f$ are H\"older continuous with exponent $0<\alpha\le1$, for all $1\le i,j\le d$.
If $\utx$, $\bU^{(j)}(\t,x)$, and $\sU^{(j)}(\t,x)$ are given by \eqref{expectn} then  the family $$\lbr u(\t,x),  \bU^{(j)}(\t,x),\sU^{(j)}(\t,x)\rbr_{j\in\Nn}$$ solves the system of linear Kuramoto-Sivashinsky-variant interacting PDEs
\begin{equation}\label{KSPDE}
 \dfrac{\partial}{\partial t_{j}}u=- \dfrac{1}{8}\Delta^2\sU^{(j)}-\dfrac12\Delta \bU^{(j)}-\dfrac{1}{2}u; \quad \t\in\intRpn,\,x\in\Rd,\ j\in\Nn, 
\end{equation}
with temporal boundary conditions given by \eqref{ksbdrycnd}.
\end{thm}
If $n=1$, then $\sU^{(j)}=\bU^{(j)}=u$, and \eqref{KSPDE} becomes the linear Kuramoto-Sivashinsky PDE (equation (3) in \cite{Aks}). 
\section{Proof of \thmref{BTBSPDE}}
\subsection{A technical lemma: the BTBS case} We start with a differentiating-under-the-integral type lemma that is needed in the proof of \thmref{BTBSPDE}.
\begin{lem}\label{diffunderint}
Let $\Wx$ be an $n$-parameter $d$-dimensional Brownian sheet starting at $x\in\Rd$ under $\P;$ and let  $f:\Rd\rightarrow\R$ be bounded and measurable such that $\pa^{2}_{kl} f$ is H\"older continuous, with exponent $0<\alpha\le1$, for $1\le k,l\le d$.  Fix $j\in\Nn$ and let
\begin{equation}
\begin{split}
\sE^{(j)}(\t,x)&\eqdef\int_{\Rpn} \left(\prod_{i\in\Nn\setminus\lbr j\rbr}s_i^{2}\right)\EP\lbk f(\Wpsx)\rbk \prod_{i=1}^{n}\KBmitzs d\s\\
\end{split}
\end{equation}
for $j\in\Nn$, $(\t,x)\in\intRpn\times\Rd$, where $\s=\lpa s_{1},\ldots,s_{n}\rpa$ and $\KBmitzs $ is the transition density of $B^{(i)}(t_{i})$.  Then  $\Delta_{x}^2\sE^{(j)}(\t,x)$ is finite and
\begin{equation}\label{diffunderint1}
\begin{split}
\Delta_{x}^2\sE^{(j)}(\t,x)&=\int_{\Rpn} \left(\prod_{i\in\Nn\setminus\lbr j\rbr}s_i^{2}\right)\lap_{x}^{2}\EP\lbk f(\Wpsx)\rbk \prod_{i=1}^{n}\KBmitzs d\s
\end{split}
\end{equation}
If we additionally assume that $\pa^{2}_{kl}f$ are bounded $($for all $1\le k,l\le d$\/$)$, then $\Delta_{x}^2\sE^{(j)}(\t,x)$ is continuous
on $\intRpn\times\Rd$.
\end{lem}
\bpf
For notational simplicity we show that 
\begin{equation}\label{toshow}
 \frac{\partial^4\sE^{(j)}}{\partial x_k^4}=\int_{\Rpn} \left(\prod_{i\in\Nn\setminus\lbr j\rbr}s_i^{2}\right)\frac{\partial^4}{\partial x_k^4}\EP\lbk f(W^{x}(\s))\rbk \prod_{i=1}^{n}\KBmitzs d\s;\  k=1,\ldots,d,\\
\end{equation}
the mixed derivatives cases follow the same steps.
In the remainder of the proof, $C$ denotes a constant that depends only on $n,d,$ and $\alpha$ and that may change its value from line to line.  Fix an arbitrary $k\in\{1,\ldots d\}$, and let 
\beq\lbl{bsden}
\bsp
\KBssxy:=\frac{1}{dy}\P[W^{x}(\s)\in dy]=\frac{\exp\lpa\dfrac{-|x-y|^2}{2\prod_{i=1}^{n}s_i}\rpa}{{\lpa 2\pi \prod_{i=1}^{n}s_i\rpa^{d/2}}} 
\end{split}
\eeq
be the transition density for the Brownian sheet $W^{x}$, starting on $\bRpn$ from $x\in\Rd$, under $\P$, and moving to $y\in\Rd$ at $\s=(s_1,\ldots,s_{n})\in\intRpn$.  Using the boundedness on $f$ and Problem 3.1 p.~254 in \cite{KS} (the case $\Rpn\times\Rd$ with $n,d>1$ is a simple extension when $f$ is bounded), the symmetry of  $\KBssxy$ in $x$ and $y$, and the facts that 
$$\lim_{y_k\to\pm\infty}f(y)\frac{\partial^3}{\partial y_k^3}\KBssxy    =\lim_{y_k\to\pm\infty}\frac{\partial}{\partial y_k}f(y)\frac{\partial^2}{\partial y_k^2}\KBssxy  =0,$$
(since $f$ is bounded and $\frac{\partial}{\partial y_k}f(y)$ is Lipschitz in $y_k$), we get
\begin{equation}
\begin{split}   \label{fstep}
&\frac{\partial^4}{\partial x_k^4}\EP f(\Wpsx)=
\left(\int_{\Rd}f(y)\frac{\partial^4}{\partial x_k^4}\KBssxy    dy\right)\\
&= \left(\int_{\Rd}f(y)\frac{\partial^4}{\partial y_k^4}\KBssxy dy\right)\\&=
\left(\int_{\Rd}\frac{\partial^2}{\partial y_k^2}f(y)\frac{\partial^2}{\partial y_k^2}\KBssxy dy\right).
\end{split}
\end{equation}
Rewriting the last term in \eqref{fstep}, and letting $h_k(y)\eqdef\partial^2 f(y)/\partial y_k^2$, we have
\begin{equation}\label{sstep}
\begin{split}
&\left(\int_{\Rd}\lpa2\pi \prod_{i=1}^{n}s_i\rpa^{-{d/2}}\left(\frac{(x_k-y_k)^2-\prod_{i=1}^{n}s_{i}}{\prod_{i=1}^{n}s_{i}^2}\right)e^{\frac{-|x-y|^2}{2\prod_{i=1}^{n}s_{i}}}h_k(y)dy\right)\\
&= \EP\left[\left(\frac{(x_k-W_k^{x_{k}}(\s))^2-\prod_{i=1}^{n}s_{i}}{\prod_{i=1}^{n}s_{i}^2}\right)h_k(\Wpsx)\right]\\
&= \EP\left[\left(\frac{(x_k-W_k^{x_{k}}(\s))^2-\prod_{i=1}^{n}s_{i}}{\prod_{i=1}^{n}s_{i}^2}\right)\left(h_k(\Wpsx)-h_k(x)\right)\right],
\end{split}
\end{equation}
where we used the fact that $\EP\left((x_k-W_k^{x_{k}}(\s))^2-\prod_{i=1}^{n}s_{i}\right)=0$ to obtain the last equality.   Now, using the Brownian sheet scaling,
we have
\begin{equation} \label{Bsc}
\bsp
&\EP\left|(x_k-W_k^{x_{k}}(\s))^2-\prod_{i=1}^{n}s_{i}\right|^2
=\prod_{i=1}^{n}s_{i}^2\EP\left|\left(\frac{W_k^{0}(\s)}{\sqrt{\prod_{i=1}^{n}s_{i}}}\right)^2-1\right|^2
\\&=\prod_{i=1}^{n}s_{i}^2\EP\left|\left({W_k^{0}(1,\ldots,1)}\right)^2-1\right|^2=C \prod_{i=1}^{n}s_{i}^2,
\end{split}
\end{equation}
for some constant $C$. Using \eqref{sstep} and Cauchy-Schwarz inequality yield
\begin{equation}
\begin{split}
 &\left|\frac{\partial^4}{\partial x_k^4}\EP f(\Wpsx) \right| \\
 &\le  \left(\EP\left(\frac{(x_k-W_k^{x_{k}}(\s))^2-\prod_{i=1}^{n}s_{i}}{\prod_{i=1}^{n}s_{i}^2}\right)^2\EP\left|h_k(\Wpsx)-h_k(x)\right|^2\right)^{1/2} \\
 &\le \frac{C}{\prod_{i=1}^{n}s_{i}}\left(\EP\left|\Wpsx-x\right|^{2\alpha}\right)^{1/2}
 = \frac{C}{\prod_{i=1}^{n}s_{i}^{1-{\alpha/2}}},
\end{split}
\label{tstep}
\end{equation}
where the next to last inequality follows from \eqref{Bsc} and the H\"{o}lder condition on $h_k$.  
Fixing an arbitrary $j\in\{1,\ldots n\}$,  we therefore have
\begin{equation}\label{diffunderintjust}
\begin{split}
&\int_{\Rpn} \left(\prod_{i\in\Nn\setminus\lbr j\rbr}s_i^{2}\right)\lab\frac{\partial^4}{\partial x_k^4}\EP\lbk f(W^{x}(\s))\rbk\rab \prod_{i=1}^{n}\KBmitzs d\s
\\&\le\frac{C}{\sqrt{2\pi t_{j}}}\int_{\Rpnmo} \left(\prod_{i\in\Nn\setminus\lbr j\rbr}s_i^{1+\frac\alpha2}\right)\prod_{i\in\Nn\setminus\{j\}}\KBmitzs \int_{\Rp}\frac{e^{-s_{j}^2/2t}}{s_{j}^{1-{\alpha/2}}}ds_{j} d\snj
\\&=\frac{C \prod_{i\in\Nn\setminus\lbr j\rbr}t_i^{\frac{2+\alpha}{4}}}{t_{j}^{\frac{2-\alpha}{4}}}<\infty;\qquad\t\in\intRpn,\ 
0<\alpha\le1.
\end{split}
\end{equation}
So, \eqref{toshow}, as well as the continuity assertion for $ {\partial^4\sE^{(j)}}/{\partial x_k^4}$, follow by a standard classical argument.
\epf

\subsection{Proof of the BTBS-PDEs system link}
Using \eqref{heatsystem} in \lemref{BSPDE} (the proof of \lemref{BSPDE} is given in the Appendix \ref{app}), \lemref{diffunderint}, and generalizing our BTP-PDEs approach in \cite{Abtp1,Abtp2}, we are now in a position to prove \thmref{BTBSPDE}. 
\bpfs{Proof of \thmref{BTBSPDE}}  It is enough to prove \thmref{BTBSPDE} (i) since a solution of \eqref{btbsystem} is automatically a solution of \eqref{btbsnonl}.  First, fix an arbitrary $j\in\{1,\ldots,n\}$ and observe that owing to the independence of $\Wx$ and $\lbr B^{(i)}\rbr_{i=1}^{n}$ as well as the independence of all the $B^{(i)}$'s  and the boundedness of $f$ together with Fubini's theorem we have
\beq\label{BTBSexp}
\bsp
u(\t,x)=\EP\lbk f\lpa\WxBnt\rpa\rbk=2^{n}\int_{\Rpn}\EP\lbk f(\Wx(\s))\rbk\prod_{i=1}^{n}\KBmitzs d\s 
\end{split}
\eeq
\beq\label{BTBSexp2}
\bsp
\sU^{(j)}(\t,x)&=\EP\lbk \lpa\prod_{i\in\Nn\setminus\lbr j\rbr}\lab B^{(i)}(t_{i})\rab\rpa^{2} f\lpa\WxBnt\rpa\rbk
\\&=2^{n}\lbk\int_{\Rpn} \lpa\prod_{i\in\Nn\setminus\lbr j\rbr}s_i\rpa^{2}\EP\lbk f\lpa\Wx_{\s}\rpa\rbk \prod_{i=1}^{n}\KBmitzs d\s\rbk
\end{split}
\eeq
for $(\t,x)\in\intRpn\times\Rd$, where $\KBmitzs $ is the transition density of $B^{(i)}(t)$.  Differentiating $\eqref{BTBSexp}$ with respect to $t_{j}$ and putting the derivative under the integral,
which is easily justified by the dominated convergence theorem, then using the fact that $\KBmjtzs  $ satisfies the heat equation
\beq\lbl{heatkernel}
\frac{\partial}{\partial t_{j}} \KBmjtzs  =\frac12 \frac{\partial^2}{\partial s_{j}^2}\KBmjtzs 
\eeq
 we have
\begin{equation}
\begin{split}
\frac{\partial}{\partial t_{j}} u(\t,x) & = 2^{n}\int_{\Rpn} \EPx\lbk f(\Ws)\rbk\frac{\partial}{\partial  t_{j}} \KBmjtzs   \prod_{i\in\Nn\setminus\lbr j\rbr}\KBmitzs d\s  \\
&=2^{n-1}\int_{\Rpn} \EPx\lbk f(\Ws)\rbk\frac{\partial^2}{\partial s_{j}^2}\KBmjtzs   \prod_{i\in\Nn\setminus\lbr j\rbr}\KBmitzs d\s 
\end{split}
\label{tder}
\end{equation}
Integrating by parts twice in $s_{j}$ and observing that the boundary terms always vanish at $\infty$ (as $s_{j}\nearrow\infty$) and that $\dstyle{(\partial/\partial s_{j}) \KBmjtzs  |_{ s=0}=0}$ but $\KBmjtzz    >0$, and then using the Brownian-sheet PDEs connection in \lemref{BSPDE} \eqref{heatsystem} successively twice we obtain
\begin{equation*}
\begin{split}
\frac{\partial}{\partial t_{j}} \utx&=2^{n-1}\int_{\Rpn} \frac{\partial^2}{\partial s_{j}^2}\EP\lbk f(\Wx_{\s}\rbk \prod_{i=1}^{n}\KBmitzs d\s 
\\&+2^{n-1}\int_{\Rpnmo}\KBmjtzz     \left.\left( {\frac {\partial }{\partial s_{j}}}\EP\lbk f(\Ws^{x})\rbk\right)\right|_{ s_{j}=0} \prod_{i\in\Nn\setminus\lbr j\rbr}\KBmitzs d\snj  \\&=2^{n-1}\int_{\Rpn} \left(\frac12\prod_{i\in\Nn\setminus\lbr j\rbr}s_i\right)^{2}\lap_{x}^{2}\EP\lbk f(\Ws^{x})\rbk \prod_{i=1}^{n}\KBmitzs d\s
\\&+\frac{2^{n-1}}{\sqrt{2\pi t_{j}}}\int_{\Rpnmo}\left(\frac12{\prod_{i\in\Nn\setminus\lbr j\rbr}}s_i\right)  \lap_{x} f(x) \prod_{i\in\Nn\setminus\lbr j\rbr}\KBmitzs d\snj  
\\&=2^{n-3}\lap_{x}^{2}\lbk\int_{\Rpn} \left(\prod_{i\in\Nn\setminus\lbr j\rbr}s_i\right)^{2}\EP\lbk f(\Ws^{x})\rbk \prod_{i=1}^{n}\KBmitzs d\s\rbk
\\&+\frac{2^{n-2}}{\sqrt{2\pi t_{j}}} \lap_{x} f(x)\int_{\Rpnmo}\left(\prod_{i\in\Nn\setminus\lbr j\rbr}s_i\right)  \prod_{i\in\Nn\setminus\lbr j\rbr}\KBmitzs d\snj  
\\&=\sqrt{\frac{\ds\prod_{i\in\Nn\setminus\lbr j\rbr}t_i}{2^{4-n}\pi^{n}t_{j}}}\lap_{x} f(x)+\frac18\lap_{x}^{2}\sU^{(j)}(\t,x),
\end{split}
\end{equation*}
where pulling  $\lap_{x}^{2}$ outside the integral in the next to last step is justified by \lemref{diffunderint} under the H\"older continuity and boundedness assumptions on the second partial derivatives of $f$.  The boundary conditions $(b)$, $(c)$, and $(d)$ follow easily from our assumptions on $f$ and from the independence of the $B^{(i)}$'s and from their independence of the Brownian sheet $\Wx$. \epfs
\section{Proof of \thmref{PDEKS}}
\subsection{A technical lemma: the Kuramoto-Sivashinsky-variant case} As in the proof of \thmref{BTBSPDE}, we need the following differentiating-under-the-integral type lemma for the proof of \thmref{PDEKS}. 
\begin{lem}\label{diffunderintks}
Assume  $f:\Rd\rightarrow\R$ is bounded and measurable such that $\pa^{2}_{kl} f$ is H\"older continuous, with exponent $0<\alpha\le1$, for $1\le k,l\le d$.  Let $v$, $\sU^{(j)}$, and $\bU^{(j)}$ be defined as in \eqref{expectn}$:$
\begin{equation}
\begin{split}
\sU^{(j)}(\t,x)&\eqdef\int_{\Rn} \left(\prod_{i\in\Nn\setminus\lbr j\rbr}s_i^{2}\right)\vsx \prod_{i=1}^{n}\KBmitzs d\s\\
\bU^{(j)}(\t,x)&\eqdef\int_{\Rn} \left(\prod_{i\in\Nn\setminus\lbr j\rbr}s_i\right)\vsx \prod_{i=1}^{n}\KBmitzs d\s
\end{split}
\end{equation}
for $j\in\Nn$, $(\t,x)\in\intRpn\times\Rd$, and $\s=\lpa s_{1},\ldots,s_{n}\rpa\in\Rn$; and where $\KBmitzs $ is the transition density of a one dimensional BM, $B^{(i)}(t_{i})$, starting at $0$.  Then  $\Delta_{x}^2\sU^{(j)}(\t,x)$ and $\Delta_{x}\bU^{(j)}(\t,x)$ are finite and
\begin{equation}\label{diffunderint2}
\begin{split}
\Delta_{x}^2\sU^{(j)}(\t,x)&\eqdef\int_{\Rn} \left(\prod_{i\in\Nn\setminus\lbr j\rbr}s_i^{2}\right)\Delta_{x}^2\vsx \prod_{i=1}^{n}\KBmitzs d\s\\
\Delta_{x}\bU^{(j)}(\t,x)&\eqdef\int_{\Rn} \left(\prod_{i\in\Nn\setminus\lbr j\rbr}s_i\right)\Delta_{x}\vsx \prod_{i=1}^{n}\KBmitzs d\s
\end{split}
\end{equation}
If we additionally assume that $\pa^{2}_{kl}f$ are bounded $($for all $1\le k,l\le d$\/$)$, then $\Delta_{x}^2\sU^{(j)}(\t,x)$ and $\Delta_{x}\bU^{(j)}(\t,x)$ are continuous on $\intRpn\times\Rd$.
\end{lem}
In light of our adaptation of the proof of Lemma 2.1 in \cite{Abtp2} to that of Lemma 2.1 in \cite{Aks} in the one parameter case ($n=1$), the proof of \lemref{diffunderintks} is a straightforward and similar adaptation of \lemref{diffunderint} above with only minor changes, and we omit these details. 
\subsection{Proof of the Kuramoto-Sivashinsky-variant interacting system link}
We now have the needed ingredients for the proof of \thmref{PDEKS}, which we present next.
\begin{pf*}{Proof of \thmref{PDEKS}}
Let $u$, $v$, $\bU^{(j)}$, and $\sU^{(j})$ be as given in \eqref{expectn}.   Differentiating $\utx$ with respect to $t$ and putting the derivative under the integral,
which is easily justified by the dominated convergence theorem, then using \eqref{heatkernel} and integrating by parts twice using the fact that the boundary terms vanish at $\pm\infty$,  we obtain
\begin{equation}\label{tder77}
\begin{split}
&\frac{\partial}{\partial t_{j}} \utx  = \int_{\Rn} \vsx\frac{\partial}{\partial  t_{j}} \KBmjtzs   \prod_{i\in\Nn\setminus\lbr j\rbr}\KBmitzs d\s 
\\&= \frac12\int_{\Rn}\vsx\frac{\partial^2}{\partial s_{j}^2}\KBmjtzs   \prod_{i\in\Nn\setminus\lbr j\rbr}\KBmitzs d\s\\
&=\frac12\int_{\Rn}\frac{\partial^2}{\partial s_{j}^2}\vsx \prod_{i\in\Nn}\KBmitzs d\s
\\&=\frac12\int_{\Rn} \prod_{i\in\Nn}\KBmitzs\left[-\frac{1}{4} \left(\prod_{i\in\Nn\setminus\lbr j\rbr}s_i\right)^{2}\Delta_{x}^2\vsx\right.
\\&\hspace{5.25cm}\left. -\left(\prod_{i\in\Nn\setminus\lbr j\rbr}s_i\right)\Delta_{x} \vsx -\vsx\right] d\s
 \\&=-\frac18\Delta_{x}^2\int_{\Rn} \prod_{i\in\Nn}\KBmitzs \left(\prod_{i\in\Nn\setminus\lbr j\rbr}s_i\right)^{2}\vsx d\s
\\&-\frac12\Delta_{x}\int_{\Rn} \prod_{i\in\Nn}\KBmitzs \left(\prod_{i\in\Nn\setminus\lbr j\rbr}s_i\right)\vsx d\s-\int_{\Rn} \prod_{i\in\Nn}\KBmitzs \vsx d\s
\\ &=-\frac18\Delta_{x}^2\sU^{(j)}(\t,x)-\frac12\Delta_{x}\bU^{(j)}(\t,x)- \frac12 u(\t,x)
\end{split}
\end{equation}
 where we have used the fact that, for each $j\in\Nn$, $v$ satisfies
\begin{equation}
\begin{split}
\frac{\partial v}{\partial s_{j}}&=\frac{\i}{2}\left(\prod_{i\in\Nn\setminus\lbr j\rbr}s_i\right)\Delta_{x} \vsx+\i\vsx\\
\frac{\partial^2 v}{\partial s_{j}^2}&=-\frac{1}{4}\left(\prod_{i\in\Nn\setminus\lbr j\rbr}s_i\right)^{2}{\Delta_{x}^2}\vsx-\left(\prod_{i\in\Nn\setminus\lbr j\rbr}s_i\right)\Delta_{x}\vsx-\vsx,
\end{split}
\label{vint}
\end{equation}
and the conditions on $f$ and \lemref{diffunderintks} to take the applications of the derivatives outside the integrals in \eqref{tder77} and \eqref{vint}. The temporal boundary conditions are trivially satisfied in light of \eqref{ksbdrycnd}.  The proof is complete.
\end{pf*}
\newpage\appendix\section{Proof of \lemref{BSPDE}}\lbl{app}
\bpfs{Existence proof}
First, let $\KBstxy$ be the transition density of an
$n$-parameter $\Rd$-valued Brownian sheet $W^{x}({\vt})=W^{x}(t_1\ldots,t_{n})$, starting on
$\bRpn$ from $x\in\Rd$, under $\P$, and moving to $y\in\Rd$ at $\t=(t_1,\ldots,t_{n})\in\intRpn$:
$$\KBstxy:=\frac{1}{dy}\P[W^{x}_{\vt}\in dy]=\frac{1}{{(2\pi t_1
\cdots t_{n})^{d/2}}} \exp\lpa\dfrac{-|x-y|^2}{2t_1\cdots t_{n}}\rpa.$$

Then, it is easy to see that for any $1\le j\le n$,
\begin{equation}
\frac{\p}{\p t_j}\KBstxy=\frac12\lap_{x} \KBstxy\prod_{i\in\Nn\setminus\lbr j\rbr}t_i.
\label{indder}
\end{equation}
Now, by the definition of $u$, the assumptions on $f$, the dominated convergence theorem,
and \eqref{indder}, we readily get:
\begin{equation*}
\begin{split}
\displaystyle\frac{\p u}{\p t_j} &=\frac{\p}{\p t_j}\int_{-\infty}^{\infty}f(y) \KBstxy dy
=\int_{-\infty}^{\infty}f(y) \frac{\p}{\p t_j}\KBstxy dy\\
&=\frac{1}{2}\lpa\prod_{i\in\Nn\setminus\lbr j\rbr}t_i\rpa\int_{-\infty}^{\infty}f(y)\lap_{x} \KBstxy dy
\\&=\frac{1}{2}\lpa\prod_{i\in\Nn\setminus\lbr j\rbr}t_i\rpa\lap_{x}\int_{-\infty}^{\infty}f(y)\KBstxy dy
=\frac{1}{2}\lpa\prod_{i\in\Nn\setminus\lbr j\rbr}t_i\rpa\lap_{x} u;
\end{split}
\end{equation*}
for $j=1,\cdots,n, \t\in\intRpn$.  By definition and by the continuity and boundedness of $f$ we easily have $\lim_{\substack{\s\to\t\\y\to x}} u(\s,y)=u(\vt,x)=f(x)$,  for $\t\in\bRpn$ and $x\in\Rd$,  proving the existence part of \lemref{BSPDE}.
\epfs
 \bpfs{Uniqueness proof}
Fix an arbitrary $j\in\{1,\ldots,n\}$ and fix arbitrary points $\t=(t_{1},\ldots,t_{j},\ldots,t_{n})\in\Rpn$ and $x\in\Rd$.   Let $\lambda$ be Lebesgue measure on $\mathscr B\lpa\Rpnmo\rpa$, let $\mathscr R=\lbr A\subset\Rpnmo;A\in\mathscr B\lpa\Rpnmo\rpa,\lambda\lpa A\rpa<\infty \rbr$.  Now, let 
\beq\lbl{wnj}
\bsp
\sWj=\lbr\sWj_{t}(A)=\lpa\sW_{t,1}^{j,x_{1}}\lpa A\rpa,\ldots, \sW_{t,d}^{j,x_{d}}\lpa A\rpa\rpa;0\le t<\infty,A\in\sR\rbr
\end{split}
\eeq
be the $d$-dimensional space-time white noise associated with an $n$-parameter $\Rd$-valued Brownian-sheet $\Wx$ on a probability space $\OFP$:
$$\sWj_{s_{j}}\lpa\ds\otimes_{i\in\Nn\setminus\{j\}}[0,t_{i}]\rpa=\Wxtnsj;\ 0\le s_{j}<\infty,$$
 (see Walsh \cite{Wa} for the case $d=1$) where $\otimes$ denotes the cartesian product of sets, and where the $j$-th parameter $s_{j}$ is treated as time and all the other parameters $\tnj=(t_{i})_{i\in\Nn\setminus\{j\}}\in\Rpnmo$ are called space (not to be confused with the spatial variable $x\in\Rd$ in the PDEs \eqref{heatsystem}, which is the starting point in $\Rd$ of $\Wx$).  Let $\lbr\sFjt\rbr_{t\ge0}$ be the filtration generated by $\sWj$; i.e., 
$$\sFjt=\sigma\lpa\lbr\sWj_{s_{j}}(A);0\le s_{j}\le t, A\in\sR\rbr\rpa;\ t\ge0.$$
As is well known, the space-time white noise $\sWj$ is a continuous orthogonal martingale measure \cite{Wa}.  In particular,  for any arbitrary fixed $\tnj\in\Rpnmo$ and $x=(x_{1},\ldots,x_{d})\in\Rd$, we have that 
$$\Wxtnsj=\sWj_{s_{j}}\lpa\ds\otimes_{i\in\Nn\setminus\{j\}}[0,t_{i}]\rpa,\sFjsj;\ 0\le s_{j}<\infty $$
is a $d$-dimensional martingale, in $s_{j},$ with coordinates quadratic covariations given by
\beq\lbl{qcov}
\lqv W^{x_{k}}_{k,\tnj}(\cdot),W^{x_{l}}_{l,\tnj}(\cdot)\rqv_{s_j}=\delta_{kl}s_j\prod_{i\in\Nn\setminus\lbr j\rbr}t_i;\ s_{j}\in\Rp,\ k,l\in\{1,\ldots,d\},
\eeq
 where $\delta_{kl}$ is the Kronecker delta ($\delta_{kl}=1$ if $k=l$ and $\delta_{kl}=0$ if $k\neq l$). 
 
 Now, fix arbitrary points $\tnj\in\Rpnmo$ and $x\in\Rd$.  Our first step is to prove
\begin{equation}\label{step1}
\begin{split}
&\Mxtnsj\eqdef\utnj(t_j-s_j,\Wxtnsj),\ \sFjsj
\mbox{ is local martingale in } s_{j} \mbox{ on } [0,t_j), 
\\&\mbox{for any fixed but arbitrary $t_{j}>0$, whenever } u\mbox{ solves part (a) in \eqref{jthheatpde}.}
\end{split}
\end{equation}
To see \eqref{step1}, we suppose that $u$ solves the $j$-th PDE in the system \eqref{heatsystem} (a) ((a) in \eqref{jthheatpde}); and we
apply It\^{o}'s formula to the function $\utnj$, in the $s_j$
variable, to get
\begin{equation}\label{whoknows1}
\begin{split}
&\utnj(t_j-s_j,\Wxotnj(s_j))-\utnj(t_j,\Wxotnj(0)) =
\int _{0}^{s_{j}}-\p_j\utnj(t_j-r_j,\Wxotnj(r_j))dr_j\\
&+\int_{0}^{s_j}\nabla\utnj(t_j-r_j,\Wxotnj(r_j))\cdot d\Wxotnj(r_j)
\\&+\frac12\left[\prod_{i\in\Nn\setminus\lbr j\rbr}t_i\right]\int_{0}^{s_j}\lap\utnj(t_j-r_j,\Wxotnj(r_j))dr_j; \mbox{ a.s. }\P,
\end{split}
\end{equation}
where $\p_{j}$ is the partial derivative with respect to the $j$-th variable, where we used the dot product notation for the second stochastic integral,  and where the third term on the right of \eqref{whoknows1} is seen by \eqref{qcov}.
This, of course, means that $\Mxtnsj\eqdef\utnj(t_j-s_j,\Wxtnsj),\ \sFjsj$ is a local martingale in $s_{j}$ on $[0,t_j)$, since the first and the third terms on the right of \eqref{whoknows1} vanish because $u$ is assumed to solve the $j$-th PDE in the system \eqref{heatsystem} (a).  So, \eqref{step1} is proved.

The last step is to show that
\begin{equation}
\mbox{If $u$ is bounded and satisfies \eqref{jthheatpde}, then
$u(\t,x)=\EP[f(\Wx(\t))]$.}
\label{step2}
\end{equation}
If $u$ is bounded and satisfies \eqref{jthheatpde} (a); then, using \eqref{step1}, $\Mtnsj$, $0\le s_j<t_j$, is a bounded martingale.  This combined with the assumption that $u$ satisfies \eqref{jthheatpde} (b) give us that the limit as $s_{j}\nearrow t_{j}$ exists and must be
$$ M^{x}_{\tnj}(t_{j}):=\lim _{s_{j}\nearrow t_{j}} \Mxtnsj=f(\Wxotnj(t_j)), \mbox{ a.s. }$$
Since $\Mxtnsj$ is uniformly integrable then 
$$\Mxtnsj=\EP\left[M^{x}_{\tnj}(t_{j})|\sFjsj\right]=\EP\left[f(\Wxotnj(t_j))|\sFjsj\right];\quad\mbox{a.s. } 0\le s_j\le t_j.$$
Taking $s_j=0$, taking expectation, and remembering the notational conventions in \notnref{papernot} we get that 
\beq\lbl{unqj}
u(\t,x)=\utnj(t_j,x)=\EP M^{x}_{\tnj}(0)=\EP\left[f(\Wxtnj)\right]=\EP\left[f(\Wx(\t))\right], 
\eeq
and the proof is complete since $j\in\{1,\ldots,n\}$, $\t=\lpa t_{1},\ldots,t_{j},\ldots,t_{n}\rpa\in\Rpn$, and $x\in\Rd$  are arbitrary.
\epfs
\section{Frequent acronyms and notations key}\lbl{B}
\begin{enumerate}\renewcommand{\labelenumi}{\Roman{enumi}.}
\item {\textbf{Acronyms}}\vspace{2mm}
\begin{enumerate}\renewcommand{\labelenumii}{(\arabic{enumii})}
\item BM: Brownian motion.
\item BTBM: Brownian-time Brownian motion. 
\item BTBS: Brownian-time Brownian sheet. 
\item BTP:  Brownian-time process.
\item KS: Kuramoto-Sivashinsky. 
\item LKSS: linear Kuramoto-Sivashinsky sheet.
\end{enumerate}
\vspace{2.5mm}
\item {\textbf{Notations}}\vspace{2mm}
\begin{enumerate}\renewcommand{\labelenumii}{(\arabic{enumii})}
\item $\Nn=\lbr 1,\ldots,n\rbr$.
\item $\KBmitzs$:  The density of a one-dimensional BM starting at $0$.
\item $\KBstxy$: The density (or kernel) of an $n$-parameter $d$-dimensional Brownian sheet.\vspace{0.5 mm}
\item $\KSstxy$: The complex kernel or density of an $n$-paramete $d$-dimensional linear Kuramoto-Sivashinsky sheet.
\end{enumerate}
\end{enumerate}

\end{document}